\documentclass[12pt,leqno]{article}
\usepackage{amsmath,amssymb,amsthm,amscd}
\numberwithin{equation}{section} \allowdisplaybreaks
\begin{document}
\newtheorem{theorem}{Theorem}[section]
\newtheorem{defin}{Definition}[section]
\newtheorem{prop}{Proposition}[section]
\newtheorem{corol}{Corollary}[section]
\newtheorem{lemma}{Lemma}[section]
\newtheorem{rem}{Remark}[section]
\newtheorem{example}{Example}[section]
\title{Tangent Dirac structures and submanifolds}
\author{{\small by}\vspace{2mm}\\Izu Vaisman}
\date{}
\maketitle
{\def\thefootnote{*}\footnotetext[1]%
{{\it 2000 Mathematics Subject Classification:} 53D17.
}}
\begin{center} \begin{minipage}{12cm}
A{\footnotesize BSTRACT. We write down the local equations that
characterize the submanifolds $N$ of a Dirac manifold $M$ which
have a normal bundle that is either a coisotropic or an isotropic
submanifold of $TM$ endowed with the tangent Dirac structure. In
the Poisson case, these formulas prove again a result of Xu: the
submanifold $N$ has a normal bundle which is a coisotropic
submanifold of $TM$ with the tangent Poisson structure iff $N$ is
a Dirac submanifold. In the presymplectic case it is the isotropy
of the normal bundle which characterizes the corresponding notion
of a Dirac submanifold. On the way, we give a simple definition of
the tangent Dirac structure, we make new remarks about it, and we
establish characteristic, local formulas for various interesting
classes of submanifolds of a Dirac manifold.}
\end{minipage}
\end{center} \vspace{5mm}
\section{Introduction} The framework of the present paper is the $C^\infty$
category. We will denote by $\Omega^k$ spaces of differential
$k$-forms, by $\chi^k$ spaces of $k$-vector fields, by $\Gamma$
spaces of differentiable cross sections of vector bundles, and we
will use the Einstein summation convention.

The Dirac structures were introduced in the study of constrained
systems and unify Poisson and presymplectic geometry \cite{C}. We
will recall their definition later. A pair $(M,D)$ that consists
of an $n$-dimensional  manifold $M$ and a Dirac structure $D$ on
$M$ is called a {\it Dirac manifold}. We started the study of
submanifolds of a Dirac manifold in \cite{V1}, where we noticed
the classes of {\it properly normalized, totally Dirac} and {\it
cosymplectic} submanifolds. These classes extend the Poisson-Dirac
with Dirac projection, Lie-Dirac (Dirac) and cosymplectic
submanifolds of a Poisson manifold, respectively
\cite{{Xu},{CF}}.

In \cite{Xu}, Xu proved that the Dirac submanifolds of the Poisson
manifold $(M,P)$ are characterized by the nice property of having
a normal bundle which is a coisotropic submanifold of the tangent
manifold $TM$ endowed with the {\it tangent Poisson structure}.
The latter is defined by the {\it complete lift} \cite{YI} of the
bivector field $P$. All the terms of Xu's result, including the
notion of a {\it tangent Dirac structure} \cite{C1}, are also
defined for Dirac manifolds, and Xu's result indicates that one
may expect interesting connections between the geometry of a
submanifold $N$ of a Dirac manifold $M$ and the geometry of a
normal bundle of $N$ in the tangent manifold $TM$. This is the
motivation of the present paper.

We will discuss the geometric configuration of Xu's result in the
general case of a Dirac manifold. The terms of the theorem are
either new or not popular, and are based on either new or not
popular geometric constructions. Accordingly, it is an objective
of the paper to explain these terms in detail. Particularly, we
recall the general construction of the vertical and complete lifts
of tensor fields from a manifold $M$ to the total space of the
tangent bundle $TM$, and the main properties of these operations
\cite{YI}. We use these lifts in order to give a simple
definition of the tangent Dirac structure and make some new
remarks about it.

Then, we turn to submanifolds. We define various classes of
submanifolds of a Dirac manifold and characterize them via local
coordinates and bases. Furthermore, we obtain the local conditions
that characterize submanifolds $N$ of $(M,D)$ with a normal bundle
$\nu N$ which is either a coisotropic or an isotropic submanifold
of $TM$. These formulas imply the result proven by Xu in the case
of Poisson manifolds. Another consequence of the established
formulas is that the analogs of Dirac submanifolds of a
presymplectic manifold $M$ are characterized by the existence of a
normal bundle which is isotropic in $TM$.
\section{Complete and vertical lifts revisited}
Let $M$ be an $n$-dimensional, differentiable manifold and $TM$ be
the total space of its tangent bundle. In the space of
differentiable functions $C^\infty(TM)$ one has the important
linear subspace $\mathbb{L}(TM)$ of the {\it fiberwise linear
functions}, the latter being functions of the form
\begin{equation}\label{fliniare} l_\alpha(x,v)=\alpha_x(v),\hspace{5mm}x\in
M, v\in T_xM, \alpha\in\Omega^1(M).\end{equation} In particular,
if we denote by $x^i$ $(i=1,...,n)$ local coordinates on $M$ and
by $v^i$ the corresponding natural coordinates on the fibers of
$TM$ (i.e., coordinates of tangent vectors with respect to the
bases $(\partial/\partial x^i)$), we have $l_{dx^i}=v^i$. Hence,
locally, $C^\infty(TM)$ is functionally spanned by the set of
functions $(\pi^*f=f\circ\pi,l_{dg})$, where $\pi:TM\rightarrow M$
is the natural projection and $f,g\in C^\infty(M)$. In what
follows the function $\pi^*f$ will be denoted again by $f$.

Two other ingredients also are of great importance in the geometry
of $TM$.

The first is the {\it Euler vector field} $E\in\chi^1(TM)$ of
infinitesimal homotheties of the fibers, which is characterized by
\begin{equation}\label{campE} Ef=0,
\;El_\alpha=l_\alpha,\hspace{5mm} f\in
C^\infty(M),\,\alpha\in\Omega^1(M), \end{equation} and has the local
expression \begin{equation}\label{campE1} E =
v^i\frac{\partial}{\partial v^i}. \end{equation}

The second is the {\it tangent structure} tensor field
$S\in\Gamma(End\,T(TM))$, which is characterized by
\begin{equation}\label{campS} (S\mathcal{X})f=0,\;
(S\mathcal{X})l_\alpha=\alpha(\pi_*\mathcal{X}),\end{equation}
where $\mathcal{X}\in\chi^1(TM),f\in
C^\infty(M),\,\alpha\in\Omega^1(M)$, and has the local expression
\begin{equation}\label{campS1} S(\xi^i\frac{\partial}{\partial x^i}+
\eta^i\frac{\partial}{\partial v^i}) = \xi^i\frac{\partial}{\partial
v^i}. \end{equation} If, as usual, we denote by
$\mathcal{V}\subseteq T(TM)$ the subbundle tangent to the fibers,
called the {\it vertical bundle}, we have
\begin{equation}\label{propS} S^2=0,\;ker\,S=im\,S=\mathcal{V},
\end{equation} and $S$ has a vanishing {\it Nijenhuis tensor}:
\begin{equation}\label{propS1}
\mathcal{N}_S(\mathcal{X},\mathcal{Y}) = [S\mathcal{X},S\mathcal{Y}]
-S[S\mathcal{X},\mathcal{Y}] - S[\mathcal{X},S\mathcal{Y}] +
S^2[\mathcal{X},\mathcal{Y}] =0. \end{equation}

Firstly, using the ingredients introduced above and denoting
$\mathcal{T}^p_q(M)=(\otimes^pTM)\otimes(\otimes^qT^*(M))$, we get
\begin{prop}\label{liftV} {\rm\cite{YI}} There exists a unique
homomorphism of real tensor algebras that sends a tensor
$t\in\mathcal{T}^p_{q,x}(M)$ to a tensor
$t^V\in\mathcal{T}^p_{q,v}(TM)$, where $v\in TM$ and $\pi(v)=x\in
M$, called the {\it vertical lift}, such that
\begin{equation}\label{liftV1}
1^V=1,\;\alpha^V=\pi^*\alpha,\;X^V=S\mathcal{X},
\end{equation} where $\alpha\in T^*_xM$, $X\in T_xM$, and
$\mathcal{X}\in T_v(TM)$ is any vector such that
$\pi_*\mathcal{X}=X$. Moreover, the vertical lift of a
differentiable tensor field is a differentiable tensor field and,
for differential forms, the vertical lift commutes with the exterior
differential. \end{prop} \begin{proof} We notice that
$\forall\lambda\in T^*_v(TM)$ there exists a unique
$\alpha_\lambda\in T^*_xM$ such that
$\alpha_\lambda(X)=\lambda(S\mathcal{X})$ $(\pi_*\mathcal{X}=X)$.
Then, we define \begin{equation}\label{defliftV} t^V(\mathcal{X}_1,
...,\mathcal{X}_q,\lambda_1,...,\lambda_p) = t(\pi_*\mathcal{X}_1,
...,\pi_*\mathcal{X}_q,\alpha_{\lambda_1},...,\alpha_{\lambda_p}).
\end{equation} The assertions about tensor fields and differential
forms follows from the fact that
in the case of a differential form $\Phi$ (functions included) one
has $\Phi^V=\pi^*\Phi$, and in the case of a vector field
$X=\xi^i(\partial/\partial x^i)$ one has
$X^V=\xi^i(\partial/\partial v^i)$. \end{proof}

Secondly, we define an operation on tensor fields known as the {\it
complete lift} \cite{YI}. For any vector field $X\in\chi^1(M)$, the
flow $exp\,tX$ lifts to a local $1$-parameter Lie group
$(exp\,tX)_*$ on the manifold $TM$, which is defined by a vector
field $X^C\in\chi^1(TM)$ called the complete lift of $X$. The local
expression of $X^C$ is \begin{equation}\label{exprXC} X^C=
\xi^i\frac{\partial}{\partial x^i} +
v^i\frac{\partial\xi^j}{\partial x^i}\frac{\partial}{\partial v^j},
\end{equation} and it is easy to check that \cite{YI}
\begin{equation}\label{caractXC}
X^Cf^V=(Xf)^V,\;X^Cl_{df}=l_{d(Xf)}, \end{equation}
\begin{equation}\label{propXC1} (X+Y)^C=X^C+Y^C,\;(fX)^C=
f^VX^C+(l_{df})X^V,\end{equation} $$\alpha^V(X^C)=(\alpha(X))^V,$$
\begin{equation}\label{propXC2} [X^C,Y^C]=[X,Y]^C,\; [X^V,Y^C] =
[X,Y]^V,\;[X^V,Y^V]=0.\end{equation}

Furthermore, one has \begin{prop}\label{propexistC} {\rm\cite{YI}}
There exists a unique homomorphism of real linear spaces which
sends a tensor field $P\in\Gamma\mathcal{T}^p_q(M)$ to a field
$P^C\in\Gamma\mathcal{T}^p_q(TM)$, called the {\it complete lift}
of $P$, such that, $\forall f\in C^\infty(M)$, $f^C=l_{df}$,
$\forall X\in\chi^1(M)$, $X^C$ is given by {\rm(\ref{caractXC})},
and
\begin{equation}\label{conddefC}
(P\otimes Q)^C=P^C\otimes Q^V+P^V\otimes
Q^C.\end{equation}\end{prop} \begin{proof} Notice that the lift
$f^C$ of a function was chosen such that for any vector field seen
as $X:M\rightarrow TM$ the pull back $f^C\circ X=Xf$. The
definition of $f^C$ and condition (\ref{conddefC}) compel us to
define the complete lift of a $1$-form $\alpha\in\Omega^1(M)$ by
\begin{equation}\label{Cpe1forme} \alpha^C(X^V)=(\alpha(X))^V,\;
\alpha^C(X^C)=(\alpha(X))^C=l_{d(\alpha(X))}.\end{equation} The
corresponding local coordinate expression is
\begin{equation}\label{Cpe1forme1} \alpha^C=
v^j\frac{\partial\alpha_i}{\partial x^j}dx^i+\alpha_idv^i,
\end{equation} and for any vector field $X:M\rightarrow TM$ the
pull back of $\alpha^C$ is given by $X^*\alpha^C=L_X\alpha$, where
$L$ denotes the Lie derivative. Finally, condition
(\ref{conddefC}) uniquely defines the complete lift of an
arbitrary tensor field because this condition is compatible with
associativity. \end{proof}

We also indicate the following properties of the complete lift:

a) \cite{YI} The complete lift of a $k$-form $\Phi\in\Omega^k(M)$
is a $k$-form $\Phi^C\in\Omega^k(TM)$ and $d\Phi^C=(d\Phi)^C$.
Indeed, a straightforward calculation shows that this condition
holds for functions and $1$-forms. Then, the condition for an
arbitrary form follows by expressing the latter locally as a real
linear combination of exterior products of $1$-forms and using
(\ref{conddefC}).

b) \cite{YI} The Lie derivative of a tensor field
$\Phi\in\Gamma\mathcal{T}^p_q(M)$ has the following lifts:
\begin{equation}\label{liftLie}
(L_X\Phi)^V=L_{X^C}\Phi^V=L_{X^V}\Phi^C,\;(L_X\Phi)^C=L_{X^C}\Phi^C,\;
L_{X^V}\Phi^V=0.
\end{equation} It is enough to check
(\ref{liftLie}) for functions, vector fields and $1$-forms, and this
can be done with the already explained formulas (\ref{liftV1}) -
(\ref{propXC2}). Then, the general result follows from
(\ref{conddefC}).

c) The complete lift of a $k$-vector field $P\in\chi^k(M)$ is a
$k$-vector field $P^C\in\chi^k(TM)$ and the Schouten-Nijenhuis
bracket satisfies the condition $[P^C,Q^C]=[P,Q]^C$. This follows
by expressing $P,Q$ locally as real linear combinations of
exterior products of vector fields, using the expression of the
Schouten-Nijenhuis bracket of such exterior products (e.g.,
\cite{V0}) and (\ref{propXC2}), (\ref{conddefC}).

As an application of property c), if $P\in\chi^2(M)$ is a Poisson
bivector field on $M$, i.e. $[P,P]=0$, then $P^C$ is a Poisson
bivector field on the manifold $TM$. The Poisson structure defined
by $P^C$ is called the {\it tangent Poisson structure} and it was
used by many authors (\cite{{C1},{W1}}, etc.).

d) The complete lifts $X^C$ of all the vector fields
$X\in\chi^1(M)$ span a generalized foliation of $TM$ such that, if
we identify $M$ with the zero section of $TM$, the leaves through
points of $M$ are the connected components of $M$ and the leaf
through $v\in TM$, $v\neq0$, is the connected component of $v$ in
$TM\backslash M$. Indeed, by (\ref{exprXC}), at a point $v=0$ we
have
$$span\{X^C\,/\,X\in\chi^1(M)\}=span\left\{\frac{\partial}{\partial x^i}=
\left(\frac{\partial}{\partial x^i}\right)^C\right\},$$ and at a point
$v=\lambda^i(\partial/\partial v^i)$ such that, for instance,
$\lambda^1\neq0$, we have
$$span\{X^C\}=span\left\{\frac{\partial}{\partial x^i}=
\left(\frac{\partial}{\partial x^i}\right)^C,
\frac{\partial}{\partial v^i}=
\left(\frac{x^1}{\lambda^1}\frac{\partial}{\partial x^i}\right)^C
-\frac{x^1}{\lambda^1}\frac{\partial}{\partial x^i}\right\}.$$

Finally, we refer the reader to \cite{GU} and the references
therein for generalizations of the lift operations discussed
above.\\
\section{Tangent Dirac structures}
Now, we will use the complete and vertical lifts in order to define
the notion of a tangent Dirac structure, first introduced by Courant
\cite{C1}, and make some new remarks about it.

The Dirac structures are defined as a class of subbundles of the
vector bundle $E(M)=TM\oplus T^*M$. The bundle $E(M)$ has several
interesting geometric objects. The first is the non degenerate
metric of zero signature
\begin{equation}\label{ginEM}
g((X,\alpha),(Y,\beta))=\frac{1}{2}(\alpha(Y)+\beta(X)),
\end{equation} where $X,Y$ are tangent vectors and $\alpha,\beta$
are tangent covectors at $x\in M$. The second is $F\in\Gamma
End(E(M))$ given by
\begin{equation}\label{FinEM} F(X,\alpha)=(X,-\alpha), \end{equation}
which is a so-called {\it para-Hermitian structure} \cite{CFG}.
The third object is the non degenerate $2$-form
\begin{equation}\label{omegainEM}
\omega((X,\alpha),(Y,\beta))=g((X,\alpha),F(Y,\beta))=
\frac{1}{2}(\alpha(Y)-\beta(X)).
\end{equation}
\begin{defin}\label{straprD} {\rm A maximal $g$-isotropic subbundle
$A\subseteq E(M)$ is called an {\it almost Dirac structure} on $M$.}
\end{defin}

The almost Dirac structure may be interpreted in terms of $TM$
alone. Namely \cite{C}, $A$ yields the pair
$(\mathcal{A},\varpi)$, where $\mathcal{A}$ is the generalized
distribution defined as the natural projection of $A$ on $TM$ and,
$\forall x\in M$, $\varpi_x\in\wedge^2\mathcal{A}^*_x$ is the
$2$-form induced by $\omega$ of (\ref{omegainEM}) ($\forall
X,Y\in\mathcal{A}_x$, the value produced by (\ref{omegainEM}) does
not depend on the choice of $\alpha,\beta$). Conversely, the pair
$(\mathcal{A},\varpi)$ allows us to reconstruct $A$ as follows
\begin{equation}\label{reconstrA}
A=\{(X,\alpha)\,/\,X\in \mathcal{A}\,\&\,\alpha|_{\mathcal{A}}=
i(X)\varpi\}.\end{equation}

The next important thing for the bundle $E(M)$ is the {\it Courant
bracket}, which is the operation defined on $\Gamma E(M)$ by
\cite{C}
\begin{equation}\label{crosetC} [(X,\alpha),(Y,\beta)] = ([X,Y],
L_X\beta-L_Y\alpha+d(\omega((X,\alpha),(Y,\beta))))\end{equation}
$$= ([X,Y],i(X)d\beta-i(Y)d\alpha +
\frac{1}{2}d(\beta(X)-\alpha(Y))),$$ where $X,Y$ are vector fields
and $\alpha,\beta$ are differential $1$-forms on $M$, $[X,Y]$ is
the usual Lie bracket and $L$ denotes the Lie derivative. The
Courant bracket is skew-symmetric but satisfies a more complicated
than the Jacobi identity.
\begin{defin}\label{strDirac} {\rm An almost Dirac structure
$D\subseteq E(M)$ is called a {\it Dirac structure} on $M$ if
$\Gamma D$ is closed by Courant brackets.}\end{defin}

In \cite{C}, it was proven that the almost Dirac structure is
Dirac iff the equivalent pair $(\mathcal{A},\varpi)$ satisfies the
following conditions: i) $\mathcal{A}$ is a generalized foliation,
ii) the form $\varpi$ is closed along the leaves of $\mathcal{A}$.
This means that the leaves of $\mathcal{A}$ are {\it presymplectic
manifolds} (of a non constant rank!) and are called the {\it
presymplectic leaves} of $D$. If the leaves are symplectic $D$ is
equivalent with a Poisson structure. Namely, if $P$ is the
corresponding Poisson bivector field, the Dirac structure is
\begin{equation}\label{strDP} D_P=\{(i(\alpha)P,\alpha)\,/\,
\alpha\in T^*M\}.\end{equation} If the leaves are the
connected components of $M$, $D$ is a presymplectic structure on
$M$ with the presymplectic form $\varpi$ such that
\begin{equation}\label{strpresD} D=D_\varpi=\{(X,i(X)\varpi)
\,/\,X\in TM\}.\end{equation}

Another fundamental property of a Dirac structure is that the
restriction of the Courant bracket (\ref{crosetC}) to $\Gamma D$
makes $D$ into a Lie algebroid where the anchor is projection on
$TM$ (e.g., \cite{V0}).\\

In \cite{C1}, a Dirac structure of $M$ was lifted to the manifold
$TM$. In what follows, we give a simpler definition of this lift.
For this purpose we look at the locally free
$\underline{C^\infty(M)}$-module sheaf $\underline{D}$ of rank $n$
of the germs of cross sections of $D$, where
$\underline{C^\infty(M)}$ is the sheaf of germs of differentiable
functions on $M$ (e.g., \cite{Ten}). Then, we get
\begin{prop}\label{costrlift} The $\underline{C^\infty(TM)}$-module
sheaf $\underline{D}^{tg}$ spanned by the germs $(X^C,\alpha^C)$,
$(X^V,\alpha^V)$, $\forall(X,\alpha)\in\underline{D}$ is locally
free of rank $2n$ and it is isomorphic with the sheaf of germs of
cross sections of a Dirac structure $D^{tg}$ on $TM$.\end{prop}
\begin{proof} Firstly, we show that the sheaf $\underline{D}^{tg}$ is
locally free of rank $2n$. If $(B_i,\epsilon_i)$ $(i=1,...,n)$ is a
local basis for the sheaf $\underline{D}$ on $M$, an arbitrary germ
$(X,\alpha)\in\underline{D}$ is of the form
$$(X,\alpha) = \sum_{i=1}^n\lambda_i(B_i,\epsilon_i),$$ whence
$$\begin{array}{l} (X^C,\alpha^C) = \sum_{i=1}^n[\lambda_i^V(B_i^C,\epsilon_i^C)
+ \lambda_i^C(B_i^V,\epsilon_i^V)],\vspace{2mm}\\ (X^V,\alpha^V) =
\sum_{i=1}^n\lambda_i^V(B_i^V,\epsilon_i^V).
\end{array}$$ This shows that
$(B_i^C,\epsilon_i^C),(B_i^V,\epsilon_i^V)$ is a local basis of
$\underline{D}^{tg}$. (If we assume that
$$B_i=b_i^j\frac{\partial}{\partial
x^j},\;\epsilon_i=\epsilon_{ij}dx^j$$  and use formulas
(\ref{exprXC}), (\ref{Cpe1forme1}), linear independence follows
from that of $(B_i,\epsilon_i)$.)

Thus \cite{Ten}, $\underline{D}^{tg}$ is isomorphic with the sheaf
of germs of cross sections of the vector bundle with local
trivialization bases $(B_i^C,\epsilon_i^C),(B_i^V,\epsilon_i^V)$,
which may be identified with a vector subbundle $D^{tg}$ of
$T(TM)$.

Furthermore, if we denote by indices $M$ and $TM$, respectively,
objects on the two manifolds, formulas (\ref{ginEM}),
(\ref{omegainEM}) and (\ref{propXC1}), (\ref{Cpe1forme})  give
\begin{equation}\label{demisotr} \begin{array}{l}
g_{TM}((X^C,\alpha^C),(Y^C,\beta^C))=(g_M((X,\alpha),(Y,\beta)))^C,
\vspace{2mm}\\

g_{TM}((X^C,\alpha^C),(Y^V,\beta^V))=(g_M((X,\alpha),(Y,\beta)))^V,
\vspace{2mm}\\

g_{TM}((X^V,\alpha^V),(Y^V,\beta^V))=0, \end{array}\end{equation}
and similar formulas relate $\omega_{TM}$ to $\omega_M$. These
formulas ensure the isotropy property for $D^{tg}$.

Finally, from (\ref{demisotr}) for $\omega$, (\ref{caractXC}) -
(\ref{propXC2}) and property b) of the complete lift given at the
end of Section 1, we get the following formulas for Courant
brackets:
\begin{equation}\label{crosetpeTM} \begin{array}{l}
[(X^C,\alpha^C),(Y^C,\beta^C)]=[(X,\alpha),(Y,\beta)]^C,\vspace{2mm}\\

[(X^C,\alpha^C),(Y^V,\beta^V)]=[(X,\alpha),(Y,\beta)]^V,\vspace{2mm}\\

[(X^V,\alpha^V),(Y^V,\beta^V)]=0.\end{array} \end{equation} These
formulas, and the fact that a Dirac structure is a Lie algebroid,
ensure that $D^{tg}$ is closed by Courant brackets. \end{proof}
\begin{defin}\label{defliftD} {\rm The Dirac structure $D^{tg}$ of
$TM$ provided by Proposition \ref{costrlift} is called the {\it
tangent Dirac structure} of the Dirac structure $D$ of
$M$.}\end{defin}

The comparison of the generating pairs of $D^{tg}$ with the bases
produced by the computations of \cite{C1} or with the alternative
definition of the reviewer of that paper \cite{deL} shows that the
tangent Dirac structure of Definition \ref{defliftD} is the same as
that of \cite{C1}.
\begin{rem}\label{obsptS} {\rm The tangent Dirac structure is
invariant by the operator $S$ of the tangent structure of the
manifold $TM$. Indeed, the action of $S$ is defined by
$S(\mathcal{X},\Xi)=(S\mathcal{X},\Xi\circ S)$ $(\mathcal{X},\Xi)\in
D^{tg})$, and the definitions of Section 1 yield $$
SX^C=X^V,\;SX^V=0,\;\alpha^C\circ S=\alpha^V,\;\alpha^V\circ S=0.$$}
\end{rem}
\begin{example}\label{liftCP} {\rm A Poisson structure of $M$ defined
by the bivector field $P$ with $[P,P]=0$ is equivalent with the
Dirac structure $D_P$ given by (\ref{strDP}). From (\ref{propXC1})
and (\ref{Cpe1forme}), it follows easily:
$$i(\alpha^V)P^C=(i(\alpha)P)^V,\;
i(\alpha^C)P^C=(i(\alpha)P)^C.$$ Hence, the tangent Dirac
structure of $D_P$ is the Poisson structure defined on $TM$ by the
bivector field $P^C$, which is the usual definition of a tangent
Poisson structure. Similarly, if $M$ has a presymplectic structure
defined by the closed $2$-form $\varpi$, this structure may be
seen as the Dirac structure $D_\varpi$ given by (\ref{strpresD})
and the tangent Dirac structure of $D_\varpi$ is the presymplectic
structure defined on $TM$ by $\varpi^C$.}
\end{example}
\begin{example}\label{strcomplgener} {\rm The construction of the
tangent Dirac structure extends to {\it complex Dirac structures}
$L\subseteq E(M)\otimes_{\mathbb{R}} \mathbb{C}$. Such a structure
$L$ is a {\it generalized complex structure} of $M$ if $L\cap \bar
L=0$ \cite{Gl}. By looking at a complex basis $(B_i,\epsilon_i)$
of $L$, it follows easily that if $L$ is a generalized complex
structure the same is true for its tangent structure $L^{tg}$.
Therefore, the tangent manifold of a generalized complex manifold
is a generalized complex manifold, again, in a canonical way. If
$M$ has a usual complex structure, $L$ is the direct sum of the
holomorphic tangent bundle of $M$ and the anti-holomorphic
cotangent bundle \cite{Gl}, and $L^{tg}$ has the similar structure
for the usual complex structure of $TM$. On the other hand, if the
generalized complex structure is defined by a symplectic form
$\omega$ on $M$, $L$ is the complexification of the Dirac
structure $D_\omega$ of (\ref{strpresD}) \cite{Gl}, and the
generalized complex structure of $TM$ is defined by the symplectic
form $\omega^C$.} \end{example}

Now, we will give some more results about the tangent Dirac
structure. \begin{prop}\label{folprestg} If $S$ is a presymplectic
leaf of $D$ on $M$ with the presymplectic form $\varpi$ and if
$v\in TS$, the presymplectic leaf of $D^{tg}$ through $v$ is the
tangent manifold $TS\subseteq TM$, and its presymplectic form is
$\varpi^C$, where the complete lift is from $S$ to $TS$.
\end{prop}
\begin{proof} Obviously, the
tangent space of the presymplectic leaf $S_{(x,v)}(D^{tg})$ of
$D^{tg}$ at a point $(x,v)\in TM$, $\pi(v)=x$, is spanned by the
vectors $X^C(x,v)$, $X^V(x,v)$, where $X$ is a vector field
tangent to the presymplectic leaf $S_x(D)$. In particular, if
$v\in TS$, we get the first part of the proposition. Furthermore,
if $\varpi$ is the presymplectic form of the leaf $S$ we have
\cite{C}
\begin{equation}\label{exprvarpi} \varpi(X,Y)=
\omega((X,\alpha),(Y,\beta)),\end{equation} where $X,Y\in TS$,
$(X,\alpha),(Y,\beta)\in D$ and the form $\omega$ is defined by
formula (\ref{omegainEM}). Then, the definitions and properties of
the complete and vertical lifts yield
\begin{equation}\label{eqpropfoi}\begin{array}{l} \varpi^C(X^C,Y^C)=
\omega_{TM}((X^C,\alpha^C),(Y^C,\beta^C)),\vspace{2mm}\\
\varpi^C(X^C,Y^V)=
\omega_{TM}((X^C,\alpha^C),(Y^V,\beta^V)),\vspace{2mm}\\
\varpi^C(X^V,Y^V)=
\omega_{TM}((X^V,\alpha^V),(Y^V,\beta^V))=0,
\end{array} \end{equation} and we are done. \end{proof}

We recall that a Poisson structure $P$ is called {\it homogeneous}
if there exists a vector field $Z$ such that
\begin{equation}\label{homogen1} L_ZP+P=0. \end{equation} It is
well known that the tangent Poisson structure $P^C$ of any Poisson
structure $P$ of $M$ is homogeneous with $Z=E$, where $E$ is the
Euler vector field (\ref{campE1}). The generalization of
homogeneity to Dirac structures $D$ is the condition
\begin{equation}\label{homogen2} (X,\alpha)\in D \Rightarrow
([Z,X]+X,L_Z\alpha)\in D, \end{equation} which reduces to
(\ref{homogen1}) in the Poisson case \cite{V1}. Now, we get
\begin{prop}\label{Dhomogen} For an arbitrary Dirac structure $D$
of $M$, the tangent Dirac structure $D^{tg}$ is homogeneous with
$Z=E$. \end{prop} \begin{proof} The Euler field $E$ satisfies
(\ref{campE}) and also
\begin{equation}\label{campE2}
[E,X^C]=0,\;[E,X^V]=-X^V,\hspace{5mm} X\in\chi^1(M),
\end{equation} which follows by easy, local coordinates
calculations. Thus, the result is proven if we show that $\forall
(X,\alpha),(Y,\beta)\in\Gamma D$ the pairs
$(X^C,L_E\alpha^C),(0,L_E\alpha^V)$ are $g_{TM}$-orthogonal to the
pairs $(Y^C,\beta^C),(Y^V,\beta^V)$. The examination of the
corresponding scalar products shows that this is the case indeed.
\end{proof}
\begin{rem}\label{pullbackD} {\rm A Dirac structure $D$ of $M$
yields a pointwisely defined pull back $\pi^*(D)$ to $TM$ by the
natural projection $\pi:TM\rightarrow M$, which is defined by
\begin{equation}\label{pisteaD}
\pi^*(D)=\{(\mathcal{X},\pi^*\alpha),\,/\,\mathcal{X}\in T(TM),
\alpha\in T^*M, \, (\pi_*\mathcal{X},\alpha)\in D\}. \end{equation}
From (\ref{pisteaD}), it follows that the equivalent, locally free
sheaf of rank $2n$ is spanned by $(X^C,\alpha^V),(\mathcal{Z},0)$
where $(X,\alpha)\in\underline{D}$ and $\mathcal{Z}$ is vertical
on $TM$. If $(B_i,\epsilon_i)$ is a local basis of
$\underline{D}$, $(B_i^C,\epsilon_i^V),(\partial/\partial v^i,0)$
is a local basis of $\underline{\pi^*(D)}$, and we see that
$\pi^*(D)$ is a differentiable Dirac structure on $TM$, which is
different from the tangent Dirac structure. In particular, if
$D=D_P$ where $P$ is a Poisson bivector field, $\pi^*(D_P)$ is not
a Poisson structure. We might say that $D^{tg}$ is the {\it
complete lift of} $D$ and $\pi^*(D)$ is the {\it vertical lift}.
The presymplectic leaves of $\pi^*(D)$ are the restrictions $T_SM$
of the tangent bundle $TM$ to the presymplectic leaves $S$ of $D$
and the corresponding presymplectic form is the vertical lift of
the presymplectic form of $S$.}\end{rem}
\section{Submanifolds of a Dirac manifold}
We begin by defining various classes of submanifolds of a Dirac
manifold. More details and motivation on that may be found in
\cite{V1}. For simplicity, all the submanifolds are assumed to be
embedded submanifolds.
\begin{defin}\label{propnorm} {\rm A submanifold $N\hookrightarrow
M$ is {\it properly normalizable} if there exists a normal bundle
$\nu N$ of $N$ such that
\begin{equation}\label{eq1propn}
(X,\alpha)\in D|_N\,\Rightarrow\, (pr_{TN}X,pr_{T^*N}\alpha)\in
D|_N,\end{equation} where the projections are defined by the
decomposition $T_NM=\nu N\oplus TN$. If (\ref{eq1propn}) holds, the
pair $(N,\nu N)$ is a {\it properly normalized submanifold} of
$M$.}\end{defin}

It follows easily that condition (\ref{eq1propn}) is equivalent with
\begin{equation}\label{eq2propn} D|_N=(D|_N\cap(TN\oplus T^*N)\oplus
D|_N\cap(\nu N\oplus\nu^*N))\end{equation} and $D|_N\cap(TN\oplus
T^*N)$ is a differentiable Dirac structure of $N$ equal to the pull
back of $D$ by the embedding of $N$ in $M$ \cite{V1}. Therefore, a
properly normalizable submanifold has a well defined induced Dirac
structure. If $D$ comes from a Poisson structure a properly
normalizable submanifold is a Poisson-Dirac submanifold which admits
a Dirac projection in the sense of \cite{CF}. If $D$ comes from a
presymplectic form $\sigma$ the submanifold $N$ is properly
normalizable iff there exists a normal bundle $\nu N$ of $N$ which
is $\sigma$-orthogonal to $N$.

In \cite{V1} we have defined an interesting invariant of a properly
normalized submanifold called the {\it second fundamental form} of
$(N,\nu N))$ in $(M,D)$. This invariant associates with every pair
$(X,\alpha),(Y,\beta)\in D|_N\cap(TN\oplus T^*N)$ a $1$-form
$B((X,\alpha),(Y,\beta))\in\nu^*N$ with the value on $Z\in\nu N$
given by
\begin{equation}\label{formadoua} B((X,\alpha),(Y,\beta))(Z) =
Z(\tilde\alpha(\tilde Y)) -\alpha([\tilde Z,\tilde Y]) +
\beta([\tilde Z,\tilde X]),\end{equation} where $\tilde Z, (\tilde
X,\tilde\alpha),(\tilde Y,\tilde\beta)$ extend
$Z,(X,\alpha),(Y,\beta)$ from $N$ to $M$ and $(\tilde
X,\tilde\alpha),(\tilde Y,\tilde\beta)$ $\in \Gamma D$. The result
is independent of the choice of the extensions because the right
hand side of (\ref{formadoua}) is $C^\infty(M)$-linear in all
arguments.

In order to define another class of submanifolds we notice the
existence of the field of subspaces along $N$
\begin{equation}\label{pseudonorm} H_x(N,D)=\{Z\in T_xM\,/\,
\exists\alpha\in ann\,T_xN\,\&\,(Z,\alpha)\in D_x\}\;(x\in N)
\end{equation} ($ann$ denotes the annihilator
space). This field may not be differentiable and the subspaces may
have various dimensions and may intersect $T_xN$. For these
reasons we will say that $H(N,D)$ is the {\it pseudo-normal field}
of $N$ with respect to $D$.
\begin{defin}\label{sbvarcosym} {\rm The submanifold
$N\hookrightarrow M$ is a {\it cosymplectic submanifold} if the
pseudo-normal field $H(N,D)$ is a differentiable normal bundle
$\nu N$ of $N$ in $M$. $H(N,D)$ will be called the {\it natural
normal bundle} of the cosymplectic submanifold $N$.}\end{defin}

In \cite{V1}, it was proven that, if $N$ is a cosymplectic
submanifold, $(N,\nu N=H(N,D))$ is a properly normalized
submanifold, that $N$ is cosymplectic iff
\begin{equation}\label{cosympl2} D|_N\cap(TN\oplus ann\,TN)=\{0\},
\end{equation} that the induced Dirac structure of a cosymplectic
submanifold is Poisson and, along $N$, one has
\begin{equation}\label{cosympl3} D|_N=\{(X,i(X)\sigma)+
(i(\lambda)\Pi,\lambda)\,/\,X\in H(N,D),\,\lambda\in T^*N\},
\end{equation} where $\Pi\in\chi^2(N)$ is the bivector field of
the induced Poisson structure and
$\sigma\in\Gamma(\wedge^2ann\,TN)$ is a $2$-form the graph of
which is $D|_N\cap(H(N,D)\oplus H^*(N,D))$. If $D$ comes from a
Poisson structure the notion of a cosymplectic submanifold is the
known one \cite{Xu}. If $D$ comes from a presymplectic form
$\sigma$ the submanifold $N$ is cosymplectic iff the pull back of
$\sigma$ to $N$ is non degenerate.

Because of (\ref{cosympl2}) we give the following definition:
\begin{defin}\label{defect} {\rm The function
$d:N\rightarrow\mathbb{Z}$ defined by $d(x)=dim[D|_N\cap(TN\oplus
ann\,TN)]$ $(x\in N)$ is called the {\it cosymplecticity default}
of the submanifold $N$ of $(M,D)$.}\end{defin}

It turns out that the second fundamental form of a cosymplectic
manifold vanishes \cite{V1}. This property is the source of
\begin{defin}\label{deftotal} {\rm A submanifold $N$ of a Dirac manifold
$(M,D)$ is a {\it totally Dirac submanifold} if it is properly
normalizable by a certain normal bundle $\nu N$ and the second
fundamental form of $(N,\nu N)$ is zero.}
\end{defin}

One can see \cite{V1} that, if $D$ comes from a Poisson structure, a
totally Dirac submanifold is just a Dirac submanifold in the sense
of Xu \cite{Xu} (in \cite{CF} these were called Lie-Dirac
submanifolds). Indeed, if $P$ is the Poisson bivector field, the
second fundamental form of $N$ becomes
$$B((i(\alpha)P,\alpha),(i(\beta)P,\beta))(Z) = (L_{\tilde
Z}P)(\alpha,\beta)$$ $\forall\alpha,\beta\in T^*N$ (then
$i(\alpha)P,i(\beta)P\in TN$ because of the proper normalization
property), and, if $B=0$, we are in the case where $P$ is {\it
soldered} to $(N,\nu N)$ \cite{V2}. Similarly, if $D=D_\sigma$ where
$\sigma$ is a presymplectic form on $M$, for the properly normalized
submanifold $(N,\nu N)$ we get
$$B((X,i(X)\sigma),(Y,i(Y)\sigma))(Z) =
(L_{\tilde Z}\sigma)(X,Y)\hspace{5mm}(X,Y\in\chi^1(N)),$$ and the
vanishing of this form, together with $\nu N\subseteq
T^{\perp_\sigma}N$ means that $\sigma$ is soldered to $(N,\nu N)$
in the sense of
\cite{V2}.\\

Finally, in accordance with the Poisson case, we also define
\begin{defin}\label{defcoiso} {\rm A submanifold $N$ of a Dirac
manifold $(M,D)$ is {\it coisotropic} if the conditions $\alpha\in
ann\,TN$ and $(X,\alpha)\in D$ imply $X\in TN$. Dually, $N$ is an
{\it isotropic submanifold} of $(M,D)$ if the conditions $X\in TN$
and $(X,\alpha)\in D$ imply $\alpha\in ann\,TN$.}\end{defin}

In the presymplectic case the coisotropy and isotropy properties
are the classical ones (i.e., $T^{\perp_\sigma}N\subseteq TN$ and
$T^{\perp_\sigma}N\supseteq TN$, respectively, where $\sigma$ is
the presymplectic form). In the Poisson case, $N$ is coisotropic
iff, $\forall x\in N$, $T_xN\cap TS$ is a coisotropic subspace of
$TS$, respectively, $N$ is an isotropic submanifold of $S$, $S$
being the symplectic leaf through $x$. Obviously, $N$ is a
coisotropic submanifold of $(M,D)$ iff its pseudo-normal field
satisfies the condition $H(N,D)\subseteq TN$.\\

Now, we shall explain how to represent a Dirac structure $D$ of
$M$ in the neighborhood of a point $x_0$ of a submanifold $N$ of
$M$ by means of local bases.

Using a tubular neighborhood of $N$ with fibers tangent to a chosen
normal bundle $\nu N$, we get local coordinates $(x^u,y^a)$
$(u=1,...,dim\,N$; $a=1,...,codim\,N)$ around $x_0$ such that $x^u$
are coordinates along $N$ and $y^a$ are coordinates along the
tubular fibers. Then the local equations of $N$ are $y^a=0$, and
\begin{equation}\label{bazexy} TN=span\left \{\frac{\partial}{\partial
x^u}\right\}_{y=0},\;\;\nu N =span\left
\{\frac{\partial}{\partial
y^a}\right\}_{y=0},\end{equation}
\begin{equation}\label{cobazexy} T^*N=span  \{dx^u\}_{y=0},\;\;
\nu^*N=span \{dy^a\}_{y=0}.\end{equation}

On the coordinate neighborhood obtained above (shrunken if
necessary), we may consider a basis of $D$ that consists of $n$
independent pairs $(B_u,\epsilon_u),(C_a,$ $\tau_a)$ where
\begin{equation}\label{bazeBC}\begin{array}{l}
B_u= b_u^s \frac{\partial}{\partial x^s} +
b^{'h}_u\frac{\partial}{\partial y^h},\; C_a= c_a^s
\frac{\partial}{\partial x^s} + c^{'h}_a\frac{\partial}{\partial
y^h},\vspace{2mm}\\ \epsilon_u=e_{us}dx^s+e'_{uh}dy^h,\;
\tau_a=t_{as}dx^s+t'_{ah}dy^h.\end{array}\end{equation} In these,
and in all the formulas that follow, our convention is that any
index of coordinates $x$ takes the same values as the index $u$
and any index of coordinates $y$ takes the same values as the
index $a$.

Of course, these bases must satisfy the conditions implied by the
definition of a Dirac structure namely, isotropy:
\begin{equation}\label{condisoBC}
\epsilon_u(B_s)+\epsilon_s(B_u)=0,\,
\epsilon_u(C_a)+\tau_a(B_u)=0,\, \tau_a(C_h)+\tau_h(C_a)=0,
\end{equation} and integrability:
\begin{equation}\label{condintegrBC}
[(B_u,\epsilon_u),(B_s,\epsilon_s)]\in D,\;
[(B_u,\epsilon_u),(C_a,\tau_a)]\in D,\end{equation}
$$[(C_a,\tau_a),(C_h,\tau_h)]\in D.$$

We also add that the tangent distribution of the presymplectic
foliation of $D$ is
\begin{equation}\label{AdeD} \mathcal{A}(D)= span\{B_u,C_a\} \end{equation}
and the presymplectic form $\varpi$ is determined by
\begin{equation}\label{varpideD} \varpi(B_u,B_s) =
\epsilon_u(B_s),\, \varpi(B_u,C_a)=\epsilon_u(C_a),\,
\varpi(C_a,C_h)=\tau_a(C_h). \end{equation}

In what follows we write down the characteristic conditions for
the various classes of submanifolds. Definition \ref{propnorm}
shows that $N$ is a properly normalizable submanifold of $(M,D)$
iff it has a normal bundle $\nu N$ such that there exist local
bases (\ref{bazeBC}) of $D$ which satisfy the conditions
\begin{equation}\label{propnormloc}
b^{'h}_u(x,0)=0,\,e'_{uh}(x,0)=0,\,c_a^s(x,0)=0,\,t_{as}(x,0)=0.
\end{equation} We call them {\it adapted bases}.

Furthermore, $N$ is cosymplectic in $(M,D)$ iff there exist $\nu
N$ and bases that satisfy (\ref{propnormloc}) and the
supplementary conditions
\begin{equation}\label{cosymplloc} e_{us}(x,0)=\delta_{us},\;c^{'h}_a(x,0)=
\delta_a^h. \end{equation} This is an immediate consequence of
formula (\ref{cosympl3}).

Back to a properly normalized submanifold $(N,\nu N)$, with
(\ref{formadoua}) we can compute the components of the second
fundamental form with respect to adapted bases, and the result is
\begin{equation}\label{forma2local}
B((B_u,\epsilon_u),(B_v,\epsilon_v))(\left.\frac{\partial}{\partial
y^a}\right|_{y=0}) =
\left(b^s_v\frac{\partial e_{us}}{\partial y^a} +
e_{vs}\frac{\partial b_u^s}{\partial y^a}
\right)_{y=0}. \end{equation} Therefore, following Definition
\ref{deftotal}, the characterization of a totally Dirac
submanifold consists of conditions (\ref{propnormloc}) and the
annulation of the components (\ref{forma2local}).
\begin{rem}\label{redemanulff2} {\rm The skew-symmetry of the
second fundamental form $B$ is an immediate consequence of the
first condition (\ref{condisoBC}). The vanishing of $B$ for a
cosymplectic submanifold follows from the first condition
(\ref{condintegrBC}). Indeed, if $N$ is a cosymplectic
submanifold, the $TN$-component of $D|_N$ behaves like a Poisson
structure (see (\ref{cosympl3})), and the fact that
\begin{equation}\label{crosetinobs}
[(B_u,\epsilon_u),(B_s,\epsilon_s)]=([B_u,B_s],L_{B_u}\epsilon_s-L_{B_s}
\epsilon_u+d(\epsilon_u(B_s)))\end{equation} belongs to $D$ along $N$ implies
the annulation of the $1$-form component of (\ref{crosetinobs})
when calculated on a vector field of the form
$$p_a^s\frac{\partial}{\partial x^s}+ q_a^h\frac{\partial}{\partial
y^h}, \hspace{5mm}p_a^u(x,0)=0,q_a^h(x,0)=\delta_a^h,$$ and
evaluated at $y=0$. This exactly yields $B=0$.}\end{rem}

Finally, Definition \ref{defcoiso} shows that a submanifold $N$ is
coisotropic if the coefficients of the formulas (\ref{bazeBC}) are
such that $\forall(\lambda^u,\lambda^{'a})$ one has
\begin{equation}\label{eqcois}
\lambda^ue_{us}(x,0)+\lambda^{'a}t_{as}(x,0)=0\,\Rightarrow\,
\lambda^ub_u^{'h}(x,0)+\lambda^{'a}c^{'h}_a(x,0)=0.\end{equation}
The condition for an isotropic submanifold is obtained by
reversing the sense of the implication in (\ref{eqcois}).\\

Now, we have all the ingredients required to discuss the result of
Xu quoted in Introduction in the framework of Dirac manifolds. We
will deduce the conditions for a submanifold $N$ of a Dirac
manifold $(M,D)$ to have a normal bundle $\nu N$ which is either a
coisotropic or an isotropic submanifold of $(TM,D^{tg})$, and
obtain some geometric conclusions of these conditions.

We consider a point $x_0\in N$, a normal bundle $\nu N$ of $N$ in
$M$, and the local coordinates and bases of formulas
(\ref{bazeBC}) around $x_0$. Then, if we denote by $(v^u,w^a)$ the
corresponding natural coordinates on the fibers of $TM$, the
submanifold $\nu N\subseteq TM$ has the local equations
$y^a=0,v^u=0$, and
\begin{equation}\label{TnuN} T(\nu
N)=span\left\{\frac{\partial}{\partial x^u},
\frac{\partial}{\partial w^a}\right\}_{y=0,v=0}, \end{equation}
\begin{equation}\label{annuN} ann(T\nu N)=
span\{dy^a,dv^u\}_{y=0,v=0}.\end{equation}

With the bases (\ref{bazeBC}) the tangent Dirac structure $D^{tg}$
is locally spanned by $(B_u^C,\epsilon_u^C),\, (C_a^C,\tau_a^C),\,
(B_u^V,\epsilon_u^V),\, (C_a^V,\tau_a^V)$, where
\begin{equation}\label{liftCBC} \begin{array}{l}
B_u^C=l_{db^s_u}\frac{\partial}{\partial v^s} +
l_{db^{'h}_u}\frac{\partial}{\partial w^h} +
b^s_u\frac{\partial}{\partial x^s} +
b_u^{'h}\frac{\partial}{\partial y^h},\vspace{2mm}\\
C_a^C=l_{dc^s_a}\frac{\partial}{\partial v^s} +
l_{dc^{'h}_a}\frac{\partial}{\partial w^h} +
c^s_a\frac{\partial}{\partial x^s} +
c_a^{'h}\frac{\partial}{\partial y^h},\vspace{2mm}\\ B_u^V
=b^s_u\frac{\partial}{\partial v^s} +
b_u^{'h}\frac{\partial}{\partial w^h},\, C_a^V
=c^s_a\frac{\partial}{\partial v^s} +
c_a^{'h}\frac{\partial}{\partial w^h}, \end{array} \end{equation}
\begin{equation}\label{liftCet} \begin{array}{l}
\epsilon^C_u= l_{de_{us}}dx^s + l_{de'_{uh}}dy^h+ e_{us}dv^s +
e'_{uh}dw^h,\vspace{2mm}\\
\tau^C_a= l_{dt_{as}}dx^s + l_{dt'_{ah}}dy^h+
t_{as}dv^s + t'_{ah}dw^h,\vspace{2mm}\\
\epsilon^V_u = e_{us}dx^s + e'_{uh}dy^h,\,
\tau^V_a= t_{as}dx^s + t'_{ah}dy^h. \end{array} \end{equation}

By writing down a linear combination of these pairs with
coefficients $\lambda^u,\nu^a,\mu^u,\xi^a$ we get a local cross
section $(\Xi,\Psi)$ of $D^{tg}$ which has the property $\Psi\in
ann(T\nu N)$ iff \begin{equation}\label{caractann}
\begin{array}{l} (\mu^ue_{us} +
\xi^at_{as})_{y=0}=0,\,
(\lambda^ue'_{uh}+\nu^at'_{ah})_{y=0}=0,\hspace{2mm}\\
\left(\lambda^u\frac{\partial e_{us}}{\partial y^h} +
\nu ^a\frac{\partial t_{as}}{\partial y^h}\right)_{y=0}=0.
\end{array}\end{equation}

The same cross section $(\Xi,\Psi)$ satisfies the condition
$\Xi\in T(\nu N)$ iff \begin{equation}\label{caracttg}
\begin{array}{l} (\mu^ub_{u}^{s} +
\xi^ac_{a}^{s})_{y=0}=0,\,
(\lambda^ub^{'h}_{u}+\nu^ac^{'h}_{a})_{y=0}=0,\hspace{2mm}\\
\left(\lambda^u\frac{\partial b^{s}_{u}}{\partial y^h} +
\nu ^a\frac{\partial c_{a}^{s}}{\partial y^h}\right)_{y=0}=0.
\end{array}\end{equation}

Therefore, we have proven
\begin{prop}\label{condgenerale} A submanifold $N$ of a Dirac
manifold $(M,D)$ has a normal bundle $\nu N$ which is coisotropic in
$(TM,D^{tg})$ iff, around the points of $N$, $D$ has bases
{\rm(\ref{bazeBC})} such that the equations {\rm(\ref{caracttg})}
are a consequence of the equations {\rm(\ref{caractann})}.
Similarly, $\nu N$ is isotropic in $(TM,D^{tg})$ iff the equations
{\rm(\ref{caractann})} are a consequence of the equations
{\rm(\ref{caracttg})}.\end{prop}

As consequences of Proposition \ref{condgenerale} we get
\begin{prop}\label{corolXu} {\rm\cite{Xu}} A submanifold $N$ of a
Poisson manifold $M$ with the Poisson bivector field $P$ is totally
Dirac iff it has a normal bundle $\nu N$ which is a coisotropic
submanifold of $(TM,P^C)$.\end{prop}
\begin{proof} Using coordinates as in (\ref{bazeBC}) we may write
\begin{equation}\label{Pincorol}
P=\frac{1}{2}P^{us}\frac{\partial}{\partial x^u} \wedge
\frac{\partial}{\partial x^s} +Q^{ua}\frac{\partial}{\partial x^u}
\wedge \frac{\partial}{\partial y^a} +
\frac{1}{2}S^{ab}\frac{\partial}{\partial y^a} \wedge
\frac{\partial}{\partial y^b}. \end{equation} Accordingly, the bases
(\ref{bazeBC}) may be taken under the form
$$B_u=P^{us}\frac{\partial}{\partial
x^s}+Q^{ua}\frac{\partial}{\partial y^a},\,\epsilon_u=dx^u,$$
$$C_a=-Q^{ua}\frac{\partial}{\partial
x^u}+S^{ab}\frac{\partial}{\partial y^b},\,\tau_a=dy^a.$$ Then,
the equations (\ref{caractann}) become $\mu^u=0,\nu^a=0$, and the
equations (\ref{caracttg}) become
$$\sum_u\mu^uP^{us}-\sum_h\xi^hQ^{sh}=0,\,
\sum_u\lambda^uQ^{ua}+\sum_h\nu^hS^{ha}=0,\,$$ $$\sum_u\lambda^u
\frac{\partial P^{us}}{\partial y^h} - \sum_a\nu^a \frac{\partial
Q^{sa}}{\partial y^h}=0,$$ for $y=0$. Obviously, the first system
implies the second iff $Q^{ua}=0, \partial P^{us}/\partial y^h$
$=0$ for $y=0$. These exactly are the conditions for $N$ to be a
(totally) Dirac submanifold \cite{Xu}.
\end{proof}

Similarly, we have \begin{prop}\label{corolpres} A submanifold $N$
of a presymplectic manifold $M$ with the closed $2$-form $\sigma$ is
totally Dirac iff it has a normal bundle $\nu N$ which is an
isotropic submanifold of $(TM,\sigma^C)$.\end{prop}
\begin{proof} With the same notation, we have
$$\sigma= \frac{1}{2}\sigma_{us}dx^u\wedge dx^s + \varphi_{ua}dx^u
\wedge dy^a + \frac{1}{2}\theta_{ah}dy^a\wedge dy^h, $$ and the
bases $$B_u=\frac{\partial}{\partial x^u},\,C_a=
\frac{\partial}{\partial
y^a},\,\epsilon_u=\sigma_{us}dx^s+\varphi_{ua}dy^a,\,\tau_a=
-\varphi_{ua}dx^u+\theta_{ah}dy^h.$$ Then, the system
(\ref{caractann}) becomes
$$\mu^u\sigma_{us}-\xi^a\varphi_{sa}=0,\,
\lambda^u\varphi_{uh}+\nu^a\theta_{ah}=0,\,
\lambda^u\frac{\partial\sigma_{us}}{\partial y^h}-
\nu^a\frac{\partial\varphi_{sa}}{\partial y^h}=0$$ for $y=0$, and
the system (\ref{caracttg}) becomes $\mu^u=0,\nu^a=0$. The latter
conditions imply the former iff $\varphi_{ua}=0$ and
$\partial\sigma_{us}/\partial y^h=0$ for $y=0$. These are the
conditions that characterize a totally Dirac submanifold of a
presymplectic manifold.\end{proof}
\begin{prop}\label{defectn}
If $N$ is a cosymplectic submanifold of the Dirac manifold
$(M,D)$, the cosymplecticity default of the natural normal bundle
$H(N,D)$ seen as a submanifold of $(TM,D^{tg})$ satisfies the
inequalities $codim\,N\leq d(x)\leq dim\,M$ $(x\in N)$.
\end{prop}
\begin{proof} The local cross sections of
$D^{tg}\cap[T(H(N,D))\oplus ann(T(H(N,$ $D)))]$ must satisfy both
(\ref{caractann}) and (\ref{caracttg}), which, modulo
(\ref{propnormloc}) and (\ref{cosymplloc}), include the conditions
$\mu^u=0,\nu^a=0$ and do not restrict the coefficients $\xi^a$.
Therefore, $D^{tg}\cap[T(H(N,D))\oplus ann(T(H(N,D)))]$ has a
basis which consists of the pairs $(C^V_a,\tau^V_a)$ and of linear
combinations of $(B^C_u,\epsilon^C_u)$.
\end{proof}
\begin{rem}\label{obsfinala} {\rm The questions discussed above may
also be considered for the vertical lift $\pi^*(D)$ defined in
Remark \ref{pullbackD}, instead of the tangent structure $D^{tg}$.
The local bases of $\pi^*(D)$ are
$(B_u^C,\epsilon_u^V),(C_a^C,\tau_a^V), (\partial/\partial v^u,0)
(\partial/\partial w^a,0)$. The conditions (\ref{caractann}) are to
be replaced by \begin{equation}\label{anninobs}
\lambda^ue_{us}(x,0)+\nu^at_{as}(x,0)=0.\end{equation} The
conditions (\ref{caracttg}) are to be replaced by
\begin{equation}\label{tginobs}
\lambda^ub^{'h}_{u}(x,0)+\nu^ac^{'h}_{a}(x,0)=0,\,\mu^u=0,\,
\left.\lambda^u\frac{\partial
b^s_u}{\partial y^h}\right|_{y=0} + \left.\nu^a\frac{\partial
c^s_a}{\partial y^h}\right|_{y=0}=0. \end{equation} From these
formulas, we see that $\nu N$ is never coisotropic in
$(TM,\pi^*(D))$.}\end{rem}
\hspace*{7.5cm}{\small \begin{tabular}{l} Department of
Mathematics\\ University of Haifa, Israel\\ E-mail:
vaisman@math.haifa.ac.il \end{tabular}}
\end{document}